\documentclass[a4paper,12pt]{amsproc}
\usepackage{amsmath,amsfonts,amssymb,amsthm,anysize}
\usepackage{enumerate}
\usepackage{epsfig}
\usepackage{supertabular}
\usepackage{graphicx}
\allowdisplaybreaks

\newcommand\Z{{\mathbb Z}}

\newcommand\F{{\mathbb F}}

\newcommand\Tr{{\mathrm{Tr}}}

\newcommand\Cay{{\mathrm{Cay}}}

\renewcommand\mod{{\mathrm{mod\, \, }}}
\newcommand\pro{{\bf Proof: }}

\theoremstyle{plain}
\newtheorem{theorem}{Theorem}[section]
\newtheorem{problem}[theorem]{Problem}
\newtheorem{lemma}[theorem]{Lemma}

\newtheorem{proposition}[theorem]{Proposition}

\numberwithin{equation}{section}

\theoremstyle{remark}

\usepackage{url, hyperref}
\hypersetup{citecolor=blue, linkcolor=blue, colorlinks=true}

\renewcommand\le{\leqslant}

\begin{document}

\title[Strongly regular graphs]{New constructions of strongly regular Cayley graphs on abelian groups}

\author{Koji Momihara}
\address{ %
Division of Natural Science\\
Faculty of Advanced Science and Technology\\
Kumamoto University\\
2-40-1 Kurokami, Kumamoto 860-8555, Japan}
\email{momihara@educ.kumamoto-u.ac.jp}
\thanks{The author acknowledges the support by 
JSPS under Grant-in-Aid for  Scientific Research (C) 20K03719.}


\begin{abstract}
Davis and Jedwab (1997) established a great construction theory unifying many previously known constructions of difference sets, relative difference sets and divisible difference sets. They introduced the concept of building blocks, which played an important role in the theory. On the other hand, Polhill (2010) gave a construction of Paley type partial difference sets (conference graphs) based on a special system of building blocks, called a covering extended building set, and proved that there exists a Paley type partial difference set in an abelian group of order $9^iv^4$ for any odd positive integer $v>1$ and any $i=0,1$. His result covers all orders of nonelementary abelian groups in which Paley type partial difference sets exist.  
In this paper, we give new constructions of strongly regular Cayley graphs on abelian groups by extending the theory of building blocks. 
The constructions are large generalizations of Polhill's construction. 
In particular, we show that for a positive integer $m$ and elementary abelian groups $G_i$, $i=1,2,\ldots,s$, of order $q_i^4$ such that $2m\,|\,q_i+1$, there exists a decomposition of the complete graph on the abelian group $G=G_1\times G_2\times \cdots\times G_s$ by strongly regular Cayley graphs with negative Latin square type parameters $(u^2,c(u+1),- u+c^2+3 c,c^2+ c)$, where $u=q_1^2q_2^2\cdots q_s^2$ and $c=(u-1)/m$. Such strongly regular decompositions were previously known only when $m=2$ or $G$ is a $p$-group.  Moreover, we find one more new infinite family of decompositions of the complete graphs by Latin square type strongly regular Cayley graphs. Thus, we 
obtain many strongly regular  graphs with new parameters. 
\end{abstract}


\maketitle

\section{Introduction}\label{sec:intro}
A {\it strongly regular graph} with parameters $(v,k,\lambda,\mu)$ is a simple and undirected graph, neither complete nor edgeless, that has the following properties:
\begin{itemize}
\item[(i)] It is a regular graph of order $v$ and valency $k$; 
\item[(ii)] For each pair of adjacent vertices $x,y$, there are $\lambda$ vertices adjacent to both $x$ and $y$; 
\item[(iii)] For each pair of nonadjacent vertices $x,y$, there are $\mu$ vertices adjacent to both $x$ and $y$.
\end{itemize}

Let $\Gamma$ be a simple and undirected graph. The adjacency matrix of $\Gamma$ is the $(0,1)$-matrix $A$ with both rows and columns indexed by the vertex set of $\Gamma$, where $A_{x,y} = 1$ when there is an edge between $x$ and $y$ in $\Gamma$ and $A_{x,y} = 0$ otherwise. A useful way to check whether a graph is strongly regular is by using the eigenvalues of its adjacency matrix. For convenience we call an eigenvalue {\it restricted} if it has an eigenvector which is not a multiple of the all-ones vector ${\bf 1}$. Note that the restricted eigenvalues of a regular connected graph with valency $k$ are the eigenvalues different from $k$. The following theorem is a well-known characterization of strongly regular graphs~\cite{BH}. 
\begin{theorem}\label{char}
For a simple graph $\Gamma$ of order $v$, neither complete nor edgeless, with adjacency matrix $A$, the following are equivalent:
\begin{enumerate}
\item $\Gamma$ is strongly regular with parameters $(v, k, \lambda, \mu)$ for certain integers $k, \lambda, \mu$,
\item $A^2 =(\lambda-\mu)A+(k-\mu) I+\mu J$ for certain real numbers $k,\lambda, \mu$, where $I, J$ are the identity matrix and the all-ones matrix, respectively, 
\item $A$ has precisely two distinct restricted eigenvalues.
\end{enumerate}
\end{theorem}
An effective method for constructing strongly regular graphs is by using Cayley graphs. Let $G$ be an additively written abelian group of order $v$, and let $D$ be a subset of $G$ such that $0_G\not\in D$ and $-D=D$, where $-D=\{-d\mid d\in D\}$. The {\it Cayley graph on $G$ with connection set $D$}, denoted ${\rm Cay}(G,D)$, is the graph $\Gamma=(G,E)$,  where $(x,y)\in E$ for $x,y\in G$ if and only if $x-y\in D$. Note that the graph is well-defined since $0_G\not\in D$ and $-D=D$. In the case when $\Cay(G,D)$ is strongly regular, the connection set $D$ is called a (regular) {\it partial difference set}. Typical examples of strongly regular Cayley graphs are the Paley graphs, the Clebsch graph, and the affine orthogonal graphs~\cite{BH}. In particular, the partial difference set 
generating a conference graph, that is, a strongly regular graph with parameters $(v,\frac{v-1}{2},\frac{v-5}{4},\frac{v-1}{4})$, is called a {\it Paley type partial difference set}. 
We now recall the following well-known lemma in algebraic graph theory~\cite{BH}. 
\begin{lemma}
Let $G$ be an abelian group and  $D$ a subset of $G$ such that $0_G\not\in D$ and $-D=D$. Then, the restricted eigenvalues of 
$\Cay(G,D)$ are given by $\psi(D):=\sum_{x\in D}\psi(x)$, $\psi\in \widehat{G}\setminus \{\psi_0\}$, where $\widehat{G}$ is the character group of $G$ and $\psi_0$ is the principal character of $G$. 
\end{lemma}
Hence, in order to see whether  $\Cay(G,D)$ is strongly regular (equivalently, $D$ is a partial difference set in $G$), we need to check that $D$ takes exactly two non-principal character values. 
The survey of Ma~\cite{Ma} contains much of what is known about partial difference sets and about connections with strongly regular graphs. 

We consider strongly regular Cayley graphs having the special parameters $(v,k,\lambda,\mu)=(u^2,c(u-\epsilon),\epsilon u+c^2-3\epsilon c,c^2-\epsilon c)$ for some integers $u,c$ and $\epsilon\in \{-1,1\}$.  These parameters are called {\it Latin square type} or {\it negative Latin square type} depending on whether  $\epsilon=1$ or $-1$. 
Recently, strongly regular Cayley graphs with Latin square type or negative Latin square type parameters have been well-studied in relation to geometric substructures, called 
{ \it$m$-ovoids} and {\it $i$-tight sets}, in finite polar spaces. See our  survey \cite{MWX2019} for recent results and about a relationship between strongly regular Cayley graphs and such geometric substructures.


In this paper, we give new constructions of strongly regular Cayley graphs on abelian groups. Our constructions are large generalizations of 
that for Paley type partial difference sets  given by Polhill~\cite{P10}. We now give a brief review of his construction. 
His construction is based on the construction theory of difference sets established by Davis and Jedwab~\cite{DJ}. 
They introduced the concept of building blocks, which played an important role in the theory.  
A {\it building block} is a subset $A$ of an abelian group $G$ satisfying that 
there exists an integer $u$ such that $|\psi(A)|\in \{0,u\}$ for any 
non-principal  character $\psi$ of $G$. Several kinds of building blocks have been studied~\cite{bjl,DO,DJ}, which have rich 
applications in Combinatorics for constructing difference sets, relative difference sets, divisible difference sets, strongly regular graphs, directed strongly regular graphs and Hadamard matrices. 
 In particular, a special system of building blocks, called a {\it covering extended building set},  
was used for constructing Menon-Hadamard difference sets~\cite{C97,WX97,X92,XC96}. Such a system is defined in an abelian  group $G$ of order $u^2$ as four subsets 
$D_i$, $i=0,1,2,3$, of $G$ such that $|D_i|=\frac{u^2-u}{2}$ for all $i$, and  $\psi(D_i)=\pm u$ for exactly one $D_i$ and $\psi(D_j)=0$ for other $D_j$ if $\psi$ is a non-principal character. The existence of such covering extended building sets in elementary abelian groups of order $u^2=q^4$ was determined by Xia~\cite{X92} for $q\equiv 3\,(\mod{4})$ and by Chen~\cite{C97} for $q\equiv 1\,(\mod{4})$. 
On the other hand, Turyn~\cite{T84}
gave a product construction for Menon-Hadamard difference sets generating regular Hadamard matrices, which can be viewed as that for covering extended building sets. 
Then, combining these results, one can claim that there exists a 
covering extended building set in a group of order $v^4$ for any odd integer $v$. On the other hand, Polhill~\cite{P10} gave a construction of Paley type 
partial difference sets from covering extended building sets. In particular, 
he proved that $((D_0\cap D_1) \cup ((G\setminus D_2)\cap (G\setminus D_3)))\setminus \{0_G\}$ forms a Paley type partial difference set under some assumptions, and obtained the following theorem. 
\begin{theorem}\label{thm:Pol}
There exists a Paley type partial difference set (and hence  a conference graph) in an abelian group of order $9^i v^4$ for any odd integer $v>1$ and any $i=0,1$. 
\end{theorem}
On the other hand, Wang~\cite{Wan1} proved that a Paley type partial difference set in a nonelementary abelian group $G$ exists only when $|G|= v^4$ or $9v^4$ for an odd integer $v>1$. Hence, Polhill's result covers all orders of nonelementary abelain groups in which 
Paley type partial difference sets exist. 
More general necessary conditions for the existence of partial difference sets were found in \cite{Wan2}. 

In this paper, we obtain the following new families of strongly regular Cayley graphs as generalizations of Polhill's result in two different ways. 
\begin{theorem}\label{thm:main1}
Let $m$ be a positive integer, and $q_i$, $i=1,2,\ldots,s$, be prime powers
such that $2m\,|\,q_i+1$. Let $G_i$, $i=1,2,\ldots,s$, be an elementary 
abelian group of order $q_i^4$. Then, there exists a strongly regular Cayley graph on $G_1\times G_2\times \cdots\times G_s$ with negative Latin square type parameters $(u^2,c(u+1),- u+c^2+3 c,c^2+ c)$, where  $u=q_1^2q_2^2\cdots q_s^2$ and $c=i(u-1)/m$ for any $i=0,1,\ldots,m$. 
\end{theorem}
\begin{theorem}\label{thm:main2}
Let $m$ be a positive integer, and $q_i$, $i=1,2,\ldots,s$, be prime powers
such that $2m\,|\,q_i+1$. Let $G_i$, $i=1,2,\ldots,s$, be an elementary 
abelian group of order $q_i^4$. Then, there exists a strongly regular Cayley graph on $G_1\times G_2\times \cdots\times G_s$ with Latin square type parameters $(u^2,c(u-1), u+c^2-3 c,c^2- c)$, where $u=q_1^2q_2^2\cdots q_s^2$ and $c$ has either of the following forms:  
\begin{itemize}
\item[1.] $c=(i+j)(u-1)/m+j$ for any $i=0,1,2,\ldots,m-2$ and $j=0,1,2$;
\item[2.] $c=(j+i)(u-1)/m+2j$ for any $i=0,1,2,\ldots,m-1$ and $j=0,1$.  
\end{itemize}
\end{theorem}
Note that infinite families of strongly regular Cayley graphs with parameters above were previously known only when $m=2$ or $G$ is a $p$-group. 
Hence, our main theorems produce many strongly regular graphs with new parameters. 

The paper is organized as follows. 
In Section~\ref{sec:semi}, we find a good system of $2m^2$ building blocks from finite fields in semi-primitive case, which is a 
generalization of  covering extended building sets found by Xia~\cite{X92}. 
Note that the case where $m=2$ corresponds to Xia's construction taking account of complements. 
In Section~\ref{sec:proper}, we give constructions of strongly regular Cayley graphs from a system of $2m^2$ building blocks satisfying special conditions defined in Section~\ref{sec:semi}, 
which are generalizations of that for Paley type partial difference sets given by Polhill~\cite{P10}. His construction was given using the notations of unions and intersections of sets. Then, the proof
became complicated even in the case where $m=2$. On the other hand,  we use group ring notations in our construction, which enables to 
simplify the proofs and generalize the argument for the case $m>2$.   
Then, in Section~\ref{sec:prod}, we give a product construction for systems of $2m^2$ building blocks satisfying those conditions, and prove our main theorems. This construction can be viewed as  a generalization of the product 
construction for Monon-Hadamard difference sets given by Turyn~\cite{T84}. 
Our generalization is completely new and nontrivial because the number of building blocks treated in Turyn's original construction is small (only eight building blocks including complements). Then, it is  difficult to find a systematic rule for arranging building blocks behind his construction.  The trick in our  generalization is to label $2m^2$ building blocks by the elements in $\Z_m^2$ in a suitable way and compose small building blocks along the operation. 
In the final section, we give remarks and an open problem. 

\section{Systems of building blocks from finite fields}\label{sec:semi}
Let $q=p^f$ be a prime power with $p$ a prime, and $\F_q$ denote the finite field of order $q$. 
Let $\psi_{\F_q}$ be the 
canonical additive character of $\F_q$, that is, the character of $(\F_q,+)$ defined by $
\psi_{\F_q}(x)=\zeta_p^{\Tr_{q/p}(x)}$, 
where $\Tr_{q/p}$ is the trace function from $\F_q$ to $\F_p$ and $\zeta_p=\exp(2\pi i/p)$. 
Then, all the characters  of 
$(\F_q,+)$ are given as $\psi_a(x)=\psi_{\F_{q}}(ax)$ for $a\in \F_q$, where 
$x\in \F_q$~\cite{LN97}.  

Let $\omega$ be a primitive element of $\F_q$, and let $N>1$ be an integer dividing $q-1$. 
Then,  $C_i^{(N,q)}=\omega^i \langle \omega^N\rangle$, $0\leq i\leq N-1$, are called the {\it $N$th cyclotomic classes} of $\F_q$. 
In particular, 
if there exists a positive integer $j$ such that $p^j\equiv -1\,(\mod{N})$, then $\Cay(\F_q, C_0^{(N,q)})$ is strongly regular.  These examples are usually called {\it semi-primitive}~\cite{BMW82,BWX}. 
\begin{lemma}{\em (\cite{BWX})} \label{lem:aaa1}
Let $p$ be a prime. Suppose that $N>2$ and there exists a positive integer 
$j$ 
such that $p^j\equiv -1\,(\mod{N})$. Choose $j$ minimum and write 
$f=2js$ for any positive integer $s$. 
Then, 
$\Cay(\F_{p^f},C_0^{(N,p^f)})$ is a strongly regular graph with parameters
\[
(v,k,\lambda,\mu)=(u^2,c(u+(-1)^{s}),-(-1)^{s} u+c^2+3(-1)^{s} c,c^2+(-1)^{s} c)
\]
with 
$u= p^{js}$ and $c=(u-(-1)^s)/N$. 
In particular, for $i=0,1,\ldots,N-1$,  
\[
\psi_{\F_{p^f}}(C_i^{(N,p^f)})=\frac{(-1)^s p^{js}-1}{N}+ 
\begin{cases}
(-1)^{s+1}p^{js},& \text{ if $\zeta_N^i=\epsilon^s$,} \\
0, & \text{ otherwise, }
\end{cases}
\]
where \[
\epsilon=\begin{cases}
-1,& \text{ if $N$ is even and $(p^j+1)/N$ is odd}, \\
1, & \text{ otherwise.}
\end{cases}
\]
\end{lemma}

In the rest of this paper, we use the following setting: let $m$ be a positive integer and $q$ be a prime power such that 
$2m$ divides $q+1$. By Lemma~\ref{lem:aaa1}, we have the following lemmas.
\begin{lemma}\label{lem:q2}
It holds that 
\[
\psi_{\F_{q^4}}(C_i^{(q^2+1,q^4)})=
\begin{cases}
q^2-1, &\mbox{ if $i=\frac{q^2+1}{2}$,}\\
-1, &\mbox{ otherwise.}
\end{cases}
\]
\end{lemma}
\pro
Apply Lemma~\ref{lem:aaa1} as $N=q^2+1$, $\epsilon=-1$ and $s=1$.  
\qed 
\begin{lemma}\label{lem:2r}
It holds that 
\[
\psi_{\F_{q^4}}(C_i^{(2m,q^4)})=
\begin{cases}
\frac{q^2-1}{2m}-q^2, &\mbox{ if $i=0$,}\\
\frac{q^2-1}{2m}, &\mbox{ otherwise.}
\end{cases}
\]
\end{lemma}
\pro
Apply Lemma~\ref{lem:aaa1} as $N=2m$ and $s$ is even.  
\qed
\vspace{0.3cm}

We now find a good system of building blocks of $\F_{q^4}$, 
which is a generalization of covering extended building sets found by Xia~\cite{X92}. Such a system will be used as starters 
of our recursion. Their properties are studied in Proposition~\ref{prop:inex}. 

We arbitrarily partition $\{2i\,|\,i=0,1,\ldots,\frac{q^2-1}{2}\}$ into one $(\frac{q^2-1}{2m}+1)$-subset and $m-1$ $(\frac{q^2-1}{2m})$-subsets, and name the $m$ subsets as $A_x$ with any labeling $x\in \Z_m$ so that $|A_0|=\frac{q^2-1}{2m}+1$. Define  $B_x=\{-a+\frac{q^2+1}{2}\,(\mod{q^2+1})\,|\,a\in A_x\}$ for $x\in \Z_m$. 
Then, by Lemma~\ref{lem:q2}, 
if $x\not=0$, 
\[
\sum_{i\in A_x}\psi_{\F_{q^4}}(\omega^j C_i^{(q^2+1,q^4)})
=
\begin{cases}
q^2-1+(-1)(|A_x|-1)(=-\frac{q^2-1}{2m}+q^2), &\mbox{ if $j\in B_x$,}\\
(-1)|A_x|(=-\frac{q^2-1}{2m}), &\mbox{ otherwise.}
\end{cases}
\]
On the other hand, if $x=0$,  
\[
\sum_{i\in A_0}\psi_{\F_{q^4}}(\omega^j C_i^{(q^2+1,q^4)} \cup \{0\})
=
\begin{cases}
q^2+(-1)(|A_0|-1)(=-\frac{q^2-1}{2m}+q^2), &\mbox{ if $j\in B_0$,}\\
(-1)|A_0|+1(=-\frac{q^2-1}{2m}), &\mbox{ otherwise.}
\end{cases}
\]
Similarly, if $x\not=0$,  
\[
\sum_{i\in B_x}\psi_{\F_{q^4}}(\omega^j C_i^{(q^2+1,q^4)})
=
\begin{cases}
-\frac{q^2-1}{2m}+q^2, &\mbox{ if $j\in A_x$,}\\
-\frac{q^2-1}{2m}, &\mbox{ otherwise.}
\end{cases}
\]
If $x=0$, 
\[
\sum_{i\in B_0}\psi_{\F_{q^4}}(\omega^j C_i^{(q^2+1,q^4)} \cup \{0\})
=
\begin{cases}
-\frac{q^2-1}{2m}+q^2, &\mbox{ if $i\in A_0$,}\\
-\frac{q^2-1}{2m}, &\mbox{ otherwise.}
\end{cases}
\]

Let $\sigma_0$ (resp. $\sigma_1$) be an arbitrarily fixed bijection from $\{2i\,|\,i=0,1,\ldots,m-1\}$ (resp. $\{2i+1\,|\,i=0,1,\ldots,m-1\}$) to 
$\Z_m$. 
For $ x, y\in \Z_m$ with $ y\not= 0$, define 
\[
\top(x, y)=C_{\sigma_1^{-1}(x)}^{(2m,q^4)}
\cup \Big(\bigcup_{i\in A_y}C_{i}^{(q^2+1,q^4)}\Big) 
\] 
and 
\[
\bot( x, y)=C_{\sigma_0^{-1}(x)}^{(2m,q^4)}
\cup \Big(\bigcup_{i\in B_y}C_{i}^{(q^2+1,q^4)}\Big). 
\]
Furthermore,  define 
\[
\top(x, 0)=C_{\sigma_1^{-1}(x)}^{(2m,q^4)}
\cup \Big(\bigcup_{i\in A_0}C_{i}^{(q^2+1,q^4)}\Big) \cup \{0\} 
\]
and
\[
\bot( x, 0)=C_{\sigma_0^{-1}( x)}^{(2m,q^4)}
\cup \Big(\bigcup_{i\in B_0}C_{i}^{(q^2+1,q^4)}\Big) \cup \{0\}. 
\]
Then, by Lemmas~\ref{lem:q2} and \ref{lem:2r},  we have 
\[
\psi_{\F_{q^4}}(\omega^i \top(x, y))=
\begin{cases}
q^2, &\mbox{ if $i\not\equiv -\sigma_1^{-1}(x)\,(\mod{2m})$ and $i\in B_y$,}\\
-q^2, &\mbox{ if $i\equiv -\sigma_1^{-1}(x)\,(\mod{2m})$ and $i\not \in B_y$,} \\
0, &\mbox{ otherwise,}
\end{cases}
\]
and 
\[
\psi_{\F_{q^4}}(\omega^i\bot(x,y))=
\begin{cases}
q^2, &\mbox{ if $i\not\equiv -\sigma_0^{-1}(x)\,(\mod{2m})$ and $i\in A_y$,}\\
-q^2, &\mbox{ if $i\equiv -\sigma_0^{-1}(x)\,(\mod{2m})$ and $i\not \in A_y$,} \\
0, &\mbox{ otherwise.}
\end{cases}
\]
Thus, $\top(x,y)$ and $\bot(x,y)$, $x,y\in \Z_m$, become 
building blocks. The system of these building blocks in the case where $m=2$ is a covering extended building set (including complements) found by Xia~\cite{X92}. 

One can directly check that the following proposition holds. 
\begin{proposition}\label{prop:inex}
Let $G=(\F_{q^4},+)$ and $u=q^2$. Then, the sets  $\top(x,y)$ and $\bot(x,y)$, where $x,y\in \Z_m$, satisfy the following conditions: 
\begin{itemize}
\item[(1)] $|\top(x,y)|=|\bot(x,y)|=\frac{u^2-u}{m}$ for any $x,y(\not= 0)\in \Z_m$, and 
$|\top(x, 0)|=|\bot(x, 0)|=\frac{u^2-u}{m}+u$ for
 any $x\in \Z_m$. 
\item[(2)] For every $\sharp\in \{\top,\bot\}$ and $x,y\in 
\Z_m$, the sets $\sharp(x-z,y+z)$, $z\in \Z_m$, partition $G$. 
\item[(3)] For any non-principal character $\psi$ of $G$, there exists a unique pair $(x',y')\in \Z_m\times \Z_m$ and $\sharp'\in \{\top,\bot\}$ such that \begin{itemize}
\item 
$\psi(\sharp'(x',y))=-u$ for any $y(\not=y')\in \Z_m$, 
\item  
$\psi(\sharp'(x,y'))=u$ for any $x(\not=x')\in \Z_m$, 
\item $\psi(\sharp'(x,y))=0$ for any other $x,y\in \Z_m$, and 
\item $\psi(\sharp(x,y))=0$ for $\sharp\in \{\top,\bot\}\setminus \{\sharp'\}$ and any $x,y\in \Z_m$.
\end{itemize}  
\item[(4)]  For each $x\in \Z_m$, there exists a subset $S_x\subseteq G\setminus \{0_G\}$ such that  
\begin{equation}\label{eq:sumS}
\sum_{y\in \Z_m}(\top(x,y)+\bot(x,y))=
m S_x+G+[0_G] 
\end{equation}
as a group ring equation in $\Z[G]$. In particular, 
$S_x=C_{\sigma_0^{-1}(x)}^{(2m,q^4)}\cup C_{\sigma_1^{-1}(x)}^{(2m,q^4)}$. 
\item[(5)] 
For each $y_1,y_2\in \Z_m$, there exists a subset $T_{y_1,y_2}\subseteq G\setminus \{0_G\}$ such that  
\begin{equation}\label{eq:sumT}
\sum_{x\in \Z_m}(\top(x,y_1)+\bot(x,y_2))=
\begin{cases}
m T_{y_1,y_2}+G-[0_G], &  \mbox{ if $y_1\not=0$ and $y_2\not=0$,}  \\
m T_{y_1,y_2}+G+(2m-1)[0_G], & \mbox{ if $y_1=y_2=0$, } \\
m T_{y_1,y_2}+G+(m-1)[0_G], &  \mbox{ otherwise, } 
\end{cases}
\end{equation}
as a group ring equation in $\Z[G]$. In particular, $T_{y_1,y_2}=\big(\bigcup_{i\in A_{y_1}}C_{i}^{(q^2+1,q^4)}\big)\cup \big(\bigcup_{i\in B_{y_2}}C_i^{(q^2+1,q^4)}\big)$. 
\end{itemize}
\end{proposition}
\section{The sets $S_x$ and $T_{y_1,y_2}$ are partial difference sets.}
\label{sec:proper}
In this section, we assume that $G$ is an abelian group of order $u^2$, not necessarily elementary abelian, containing $2m^2$ building blocks $\top(x,y)$ and $\bot(x,y)$, where $x,y\in \Z_m$, satisfying the conditions~(1), (2), (3) and either (4) or (5) in Proposition~\ref{prop:inex}. (Regard $G$ in Proposition~\ref{prop:inex} as any abelian group of order $u^2$.) Then, we prove that the sets $S_x$ and $T_{y_1,y_2}$ defined in \eqref{eq:sumS} and \eqref{eq:sumT} are actually partial difference sets in $G$. These constructions can be viewed as 
generalizations of that for Paley type partial difference sets given by Polhill~\cite{P10}. In particular, he treated only $S_x$ in the case where $m=2$. 
\begin{proposition}\label{prop:1234}
Let $G$ be an abelian group of order $u^2$ containing $2m^2$  building blocks $\top(x,y)$ and $\bot(x,y)$, where $x,y\in \Z_m$, satisfying 
the conditions (1), (2), (3) and (4). Then, the sets $S_x$, $x\in \Z_m$, defined in \eqref{eq:sumS} satisfy the following: 
\begin{itemize}
\item[(i)] $|S_x|=\frac{u^2-1}{m}$ for any $x\in \Z_m$; 
\item[(ii)] $S_x$, $x\in \Z_m$, are pairwise disjoint;
\item[(iii)] $S_x$ takes exactly two non-principal character values 
$\frac{u-1}{m}$ and $\frac{-(m-1)u-1}{m}$.  
\end{itemize}
Hence, each $S_x$ forms a partial difference set in $G$. 
\end{proposition}
\pro 
(i) Since $\sum_{y\in \Z_m}(|\top(x,y)|+|\bot(x,y)|)=2u^2$ by the condition~(1), we have $|S_x|=\frac{u^2-1}{m}$. 

(ii) By the condition~(2), $\sum_{x\in \Z_m}\sum_{y\in \Z_m}(\top(x,y)+\bot(x,y))=2mG$. Then,   
the condition~(4) implies that  $S_x$, $x\in \Z_m$, are pairwise disjoint. 

(iii) Let $\psi$ be a non-principal character of $G$. Let $(x',y') \in \Z_m\times \Z_m$ and $\sharp'\in \{\top,\bot\}$ satisfy the 
condition~(3). If $x=x'$, 
\begin{align*}
\psi(S_x)=&\,\frac{1}{m}(-\psi(G)-\psi(0_G)+
\sum_{y\in \Z_m}
(\psi(\top(x,y))+\psi(\bot(x,y))))\\
=&\,\frac{1}{m}(-1+\sum_{y(\not=y')\in \Z_m}\psi(\sharp'(x',y)))=\frac{-1-(m-1)u}{m}. 
\end{align*}
If $x\not=x'$, 
\begin{align*}
\psi(S_x)=&\,\frac{1}{m}(-\psi(G)-\psi(0_G)+
\sum_{y\in \Z_m}
(\psi(\top(x,y))+\psi(\bot(x,y))))\\
=&\,\frac{1}{m}(-1+\psi(\sharp'(x,y')))=\frac{-1+u}{m}. 
\end{align*}
This completes the proof of the proposition. 
\qed


\begin{proposition}\label{prop:1235}
Let $G$ be an abelian group of order $u^2$ containing $2m^2$  building blocks  $\top(x,y)$ and $\bot(x,y)$, where $x,y\in \Z_m$, satisfying 
the conditions (1), (2), (3) and (5). Then, the sets $T_{y_1,y_2}$, $y_1,y_2\in \Z_m$, defined in \eqref{eq:sumT} satisfy the following: 
\begin{itemize}
\item[(i)] 
\[
|T_{y_1,y_2}|=
\begin{cases}
\frac{(u-1)^2}{m},&\mbox{ if $y_1\not=0$ and $y_2\not=0$,}\\
\frac{(u-1)^2}{m}+2(u-1),& \mbox{ if $y_1=y_2=0$, }\\
\frac{(u-1)^2}{m}+u-1,&\mbox{ otherwise.}
\end{cases}
\] 
\item[(ii)] $T_{y_1,y_1+y_2}$, $y_1\in \Z_m$, are pairwise disjoint;
\item[(iii)] $T_{y_1,y_2}$ takes exactly two non-principal character values 
$\frac{\delta_{y_1,y_2}-u}{m}$ and $\frac{\delta_{y_1,y_2}+(m-1)u}{m}$, 
where \[
\delta_{y_1,y_2}=
\begin{cases}
1,&\mbox{ if $y_1\not=0$ and $y_2\not=0$,}\\
1-2m,& \mbox{ if $y_1=y_2=0$,}\\
1-m,&\mbox{ otherwise.}
\end{cases}
\]
\end{itemize}
Hence, each $T_{y_1,y_2}$ forms a partial difference set in $G$.  
\end{proposition}
\pro 
(i) By the condition~(1),  
\[
\sum_{x\in \Z_m}(|\top(x,y_1)|+|\bot(x,y_2)|)=
\begin{cases}
2u^2-2u,&\mbox{ if $y_1\not=0$ and $y_2\not=0$,}\\
2u^2+(2m-2)u, & \mbox{ if $y_1=y_2=0$,}\\
2u^2+(m-2)u,&\mbox{ otherwise.}
\end{cases}
\] 
Then, by the condition~(5), the conclusion follows. 

(ii) By the condition~(2),  $\sum_{x\in \Z_m}\sum_{y_1\in \Z_m}(\top(x,y_1)+\bot(x,y_1+y_2))=2mG$. Then, 
the condition~(5) implies that the sets $T_{y_1,y_1+y_2}$, $y_1\in \Z_m$, are pairwise disjoint. 

(iii) Let $\psi$ be a non-principal character of $G$. Let $(x',y') \in \Z_m\times \Z_m$ and $\sharp'\in \{\top,\bot\}$ satisfy the condition~(3). 
If $\sharp'=\top$ and $y_1=y'$ or  $\sharp'=\bot$ and $y_2=y'$, then 
\begin{align*}
\psi(T_{y_1,y_2})=&\,\frac{1}{m}(-\psi(G)+\delta_{y_1,y_2}\psi(0_G)+
\sum_{x\in \Z_m}
(\psi(\top(x,y_1))+\psi(\bot(x,y_2))))\\
=&\,\frac{1}{m}(\delta_{y_1,y_2}+\sum_{x(\not=x')\in \Z_m}\psi(\sharp'(x,y')))=\frac{\delta_{y_1,y_2}+(m-1)u}{m}. 
\end{align*}
Otherwise, 
\begin{align*}
\psi(T_{y_1,y_2})=
\frac{1}{m}(\delta_{y_1,y_2}+\psi(\sharp'(x',y)))=\frac{\delta_{y_1,y_2}-u}{m}, 
\end{align*}
where $y=y_1$ or $y_2$ with $y\not=y'$.  
This completes the proof of the proposition.  
\qed

\section{Product construction for $\top(x,y)$ and $\bot(x,y)$}\label{sec:prod}
In this section, we give a product construction for $2m^2$ building blocks $\top(x,y)$ and $\bot(x,y)$, $x,y\in \Z_m$, satisfying conditions (1), (2), (3) and either (4) or (5) in  an abelian group $G$.  

Let $G_i$, $i=1,2$, be abelian groups of order $u^2$ and $v^2$, respectively. Assume that each $G_i$ contains building blocks $\top_i(x,y)$ and $\bot_i(x,y)$, where $x,y\in \Z_m$, satisfying the conditions (1), (2), (3)  and either (4) or (5). 
Define subsets $\top(x,y)$ and $\bot(x,y)$ in $G_1\times G_2$ as follows:  
\begin{align}
\top(x,y)=&\,\sum_{a, b\in \Z_m}
\top_1(x-a,y- b)\times 
(\bot_2(x+y- a, a)\cap \top_2(x+y- b, b)),\label{eq:aaa1}\\
\bot(x,y)=&\,\sum_{a,b\in \Z_m}
\bot_1(x-b,y- a)\times 
(\bot_2(x+y- a, a)\cap \top_2(x+y- b, b)), 
\label{eq:aaa2}
\end{align}
where $x,y\in \Z_m$. We will prove that $\top(x,y)$ and $\bot(x,y)$ defined above satisfy the conditions (1), (2), (3)  and either (4) or (5). This is a generalization of Turyn's product construction~\cite{T84} for covering extended building sets. In particular, the construction in the case where $m=2$ corresponds to Turyn's construction. 
\subsection{The conditions (1), (2) and (3)}
In this subsection, we show that $\top(x,y)$ and $\bot(x,y)$, $x,y\in \Z_m$, defined in \eqref{eq:aaa1} and \eqref{eq:aaa2} satisfy the conditions (1), (2) and (3). 
\begin{lemma}\label{lem:11}
The subsets $\top(x,y)$ and $\bot(x,y)$, where $x,y\in \Z_m$, defined in \eqref{eq:aaa1} and \eqref{eq:aaa2} satisfy the condition~(1). 
\end{lemma}
\pro 
By the condition~(1) for $G_1$, 
we have 
\begin{align*}
|\top(x,y)|=&\Big(\frac{u^2-u}{m}\Big)\sum_{a\in \Z_m}\sum_{b(\not=y)\in \Z_m}
|\bot_2(x+y-a,a)\cap \top_2(x+y-b,b)|\\
&\, +\Big(\frac{u^2-u}{m}+u\Big)\sum_{a\in \Z_m}
|\bot_2(x+y-a,a)\cap \top_2(x,y)|\\
=&\Big(\frac{u^2-u}{m}\Big)\sum_{a,b\in \Z_m}
|\bot_2(x+y-a,a)\cap \top_2(x+y-b,b)|\\
&\, +u\sum_{a\in \Z_m}
|\bot_2(x+y-a,a)\cap \top_2(x,y)|. 
\end{align*}
Then, by the conditions~(1) and (2) for $G_2$, we have 
\begin{align*}
|\top(x,y)|
=&\Big(\frac{u^2-u}{m}\Big)v^2
+u|\top_2(x,y))|=
\begin{cases}
\frac{u^2v^2-uv}{m}+uv,& \mbox{ if $y=0$,}\\
\frac{u^2v^2-uv}{m},& \mbox{ otherwise}.
\end{cases} 
\end{align*}
Similarly, we have 
\begin{align*}
|\bot(x,y)|=\Big(\frac{u^2-u}{m}\Big)v^2
+u|\bot_2(x,y))|
=\begin{cases}
\frac{u^2v^2-uv}{m}+uv,& \mbox{ if $y=0$,}\\
\frac{u^2v^2-uv}{m},& \mbox{ otherwise}.
\end{cases} 
\end{align*}
This completes the proof. \qed

\begin{lemma}\label{lem:21}
The subsets $\top(x,y)$ and $\bot(x,y)$, where $x,y\in \Z_m$, defined in \eqref{eq:aaa1} and \eqref{eq:aaa2} satisfy the condition~(2). 
\end{lemma}
\pro 
Let $x,y\in \Z_m$. Then, 
\begin{align*}
&\sum_{ s\in \Z_m}\top(x- s,y+ s)\\
=&\,\sum_{ s,a,b\in \Z_m}
\top_1(x-a- s,y-b+ s)
\times 
(\bot_2(x+y-a,a)\cap \top_2(x+y-b,b)).
\end{align*}
Since $\sum_{s\in \Z_m}\top_1(x-a- s,y-b+ s)=G_1$ and  $\sum_{a\in \Z_m}\bot_2(x+y-a,a)=\sum_{b\in \Z_m}\top_2(x+y-b,b)=G_2$, we have 
\begin{align*}
\sum_{s\in \Z_m}\top(x-s,y+ s)=&\,\sum_{a,b\in \Z_m}
G_1
\times (\bot_2(x+y-a,a)\cap \top_2(x+y-b,b))\\
=&\,\sum_{b\in \Z_m}
G_1
\times 
(G_2\cap \top_2(x+y-b,b))=G_1\times G_2. 
\end{align*}
Similarly, we have $\sum_{s\in \Z_m}\bot(x-s,y+ s)=G_1\times G_2$. 
This completes the proof of the lemma. 
\qed

\begin{lemma}\label{lem:31}
The subsets $\top(x,y)$ and $\bot(x,y)$, where $x,y\in \Z_m$, defined in \eqref{eq:aaa1} and \eqref{eq:aaa2} satisfy the condition~(3).
\end{lemma} 
\pro  
Let $\psi$ be a non-principal character of $G_1\times G_2$. 
We can assume that 
$\psi=(\psi_1,\psi_2) \in G_1^\perp\times G_2^\perp$. 
We compute the character values of $\top(x,y)$ and $\bot(x,y)$ dividing into two cases: (i) $\psi_1$ is principal and $\psi_2$ is non-principal; and (ii) $\psi_1$ is non-principal. 

(i) We assume that $\psi_1$ is principal and $\psi_2$ is non-principal. Then, 
\begin{align*}
\psi(\top(x,y))&=\Big(\frac{u^2-u}{m}\Big)\sum_{a\in \Z_m}
\sum_{b(\not=y)\in \Z_m} 
\psi_2(\bot_2(x+y-a,a)\cap \top_2(x+y-b,b))\\
&\hspace{2cm}+\Big(\frac{u^2-u}{m}+u\Big)\sum_{a\in \Z_m}
\psi_2(\bot_2(x+y-a,a)\cap \top_2(x,y))\\
&=\Big(\frac{u^2-u}{m}\Big)\sum_{a,b\in \Z_m} 
\psi_2(\bot_2(x+y-a,a)\cap \top_2(x+y-b,b))\\
&\hspace{2cm}+u\sum_{a\in \Z_m}
\psi_2(\bot_2(x+y-a,a)\cap \top_2(x,y))\\
&=\Big(\frac{u^2-u}{m}\Big)\psi_2(G_2)+u
\psi_2(\top_2(x,y))=u\psi_2(\top_2(x,y))\in \{0,uv,-uv\}. 
\end{align*}
Similarly, we have 
\begin{align*}
\psi(\bot(x,y))
=u\psi_2(\bot_2(x,y))\in \{0,uv,-uv\}. 
\end{align*}

(ii) Next, we assume that $\psi_1$ is non-principal. 
Let $x',y'\in \Z_m$ and $\sharp\in \{\top_1,\bot_1\}$ satisfy the condition~(3) for $G_1$. If $\sharp=\bot_1$, $\psi(\top(x,y))=0$ since  $\psi_1(\top_1(x-a,y-b))=0$ for any 
$a,b\in \Z_m$. If $\sharp=\top_1$,  we have 
\begin{align*}
\psi(\top(x,y))&=
\sum_{a,b\in \Z_m}
\psi_1(\top_1(x-a,y-b))\times 
\psi_2(\bot_2(x+y-a,a)\cap \top_2(x+y-b,b))\\
&=-u\sum_{b(\not=y-y')\in \Z_m}\psi_2(\bot_2(x'+y,x-x')\cap \top_2(x+y-b,b))\\
&\hspace{1.3cm}+u\sum_{a(\not=x- x')\in \Z_m}\psi_2(\bot_2(x+y-a,a)\cap \top_2(x+y',y-y'))\\
&=-u\psi_2(\bot_2(x'+y,x-x'))+u\psi_2(\top_2(x+y',y-y')). 
\end{align*}
Here, if $\psi_2$ is principal, 
\[
\psi(\top(x,y))=uv(-\epsilon_{x',x}+\epsilon_{y',y}) \in \{0,uv,-uv\}, 
\]
where $\epsilon_{a,b}=1$ or $0$ depending on whether $a=b$ or not. 
We now assume that $\psi_2$ is non-principal.  
Let $x'',y''\in \Z_m$ and $\sharp'\in \{\top_2,\bot_2\}$ satisfy the condition~(3) for $G_2$. If $\sharp'=\bot_2$, 
\begin{align*}
\psi(\top(x,y))=&
-u\psi_2(\bot_2(x'+y,x-x'))\\
=&\begin{cases}
uv,& \mbox{if $x'+y=x''$ and $x-x'\not=y''$,}\\
-uv,& \mbox{if $x'+y\not=x''$ and $x-x'=y''$,}\\
0,& \mbox{ otherwise.} 
\end{cases}
\end{align*}
If $\sharp'=\top_2$, 
\begin{align*}
\psi(\top(x,y))=&
u\psi_2(\top_2(x+y',y-y'))\\
=&\begin{cases}
-uv,& \mbox{if $x+y'=x''$ and $y-y'\not=y''$,}\\
uv,& \mbox{if $x+y'\not=x''$ and $y-y'=y''$,}\\
0,& \mbox{ otherwise.} 
\end{cases}
\end{align*}

Similarly,  if $\sharp=\top_1$, $\psi(\bot(x,y))=0$ since  $\psi_1(\bot_1(x-b,y-a))=0$ for any 
$a,b\in \Z_m$. Assume that $\sharp=\bot_1$. 
Then, 
\begin{align*}
\psi(\bot(x,y))
=-u\psi_2(\top_2(x'+y,x-x'))+u\psi_2(\bot_2(x+y',y-y')). 
\end{align*}
If $\psi_2$ is principal, 
\begin{align*}
\psi(\bot(x,y))
=
uv(-\epsilon_{x',x}+\epsilon_{y',y}) \in \{0,uv,-uv\}. 
\end{align*}
We next assume that $\psi_2$ is non-principal.  
Let $x'',y''\in \Z_m$ and $\sharp'\in \{\top_2,\bot_2\}$ satisfy   the condition~(3)  for $G_2$. If $\sharp'=\bot_2$, 
\begin{align*}
\psi(\bot(x,y))=&
u\psi_2(\bot_2(x+y',y-y'))\\
=&\begin{cases}
-uv,& \mbox{if $x+y'=x''$ and $y-y'\not=y''$,}\\
uv,& \mbox{if $x+y'\not=x''$ and $y-y'=y''$,}\\
0,& \mbox{ otherwise.} 
\end{cases}
\end{align*}
If $\sharp'=\top_2$, 
\begin{align*}
\psi(\bot(x,y))=&
-u\psi_2(\top_2(x'+y,x-x'))\\
=&
\begin{cases}
uv,& \mbox{if $x'+y=x''$ and $x-x'\not=y''$,}\\
-uv,& \mbox{if $x'+y\not=x''$ and $x-x'=y''$,}\\
0,& \mbox{ otherwise.} 
\end{cases}
\end{align*}

Then, 
it is straightforward to check that the sets $\top(x,y)$ and $\bot(x,y)$, where $x,y\in \Z_m$,  
satisfy the condition~(3).
\qed 

\begin{theorem}\label{thm:1half}
Let $m$ be a positive integer and $q_i$, $i=1,2,\ldots,s$, be prime powers 
such that $2m\,|\,q_i+1$. Let $G_i=(\F_{q_i^4},+)$, $1\le i\le s$, and $G=
G_1\times G_2\times \cdots \times G_s$. Then, there are $2m^2$ building  blocks $\top(x,y)$ and $\bot(x,y)$, $x,y\in \Z_m$, of $G$  satisfying the conditions~(1), (2) and (3). 
\end{theorem}
\pro
By Proposition~\ref{prop:inex}, each $G_i$ contains $2m^2$ building  blocks  satisfying the conditions~(1), (2) and (3). Apply our product construction recursively  to those examples in $G_i$ for $i=1,2,\ldots,s$. Then, by Lemmas~\ref{lem:11}, \ref{lem:21} and 
\ref{lem:31}, the conclusion follows. 
\qed
%
\subsection{The conditions (4) and (5)}
In this subsection, we show that $\top(x,y)$ and $\bot(x,y)$, $x,y\in \Z_m$, defined in \eqref{eq:aaa1} and \eqref{eq:aaa2} satisfy the conditions (4) and (5). It seems to be complicated to prove them directly. 
Therefore, throughout this subsection, in addition to the assumption in the previous subsection, 
we restrict $G_2=(\F_{q^4},+)$ with $2m\,|\,q+1$, and $\top_2(x,y)$ and $\bot_2(x,y)$, $x,y\in \Z_m$, are defined as in 
Section~\ref{sec:semi}.  

In the following lemma, we assume that $\top_1(x,y)$ and $\bot_1(x,y)$, $x,y\in \Z_m$,
satisfy the conditions~(1), (2), (3) and (4). 
\begin{lemma}\label{lem:41}
The subsets $\top(x,y)$ and $\bot(x,y)$, where $x,y\in \Z_m$, defined in \eqref{eq:aaa1} and \eqref{eq:aaa2} satisfy the condition~(4). 
\end{lemma}
\pro 
Let $(a,b)\in G_1\times G_2$, where $a$ is arbitrarily chosen. 

First, we assume that $b\in C_s^{(2m,q^4)}\cap C_t^{(q^2+1,q^4)}$ with $s$ and $t$ even. Then, there exists $(x_0,y_0)\in \Z_m\times \Z_m$ uniquely  such that  
$b\in \bot_2(x_0,y)\cap \top_{2}(x,y_0)$ for all $x,y\in \Z_m$. 
We now compute the coefficient of  $(a,b)$ in $\sum_{y\in \Z_m}(\top(x,y)+\bot(x,y))$. The $(a,b)$ may appear in the partial summation 
\begin{align*}
\sum_{y\in \Z_m}(\top_1(x_0-y,y-y_0)+\bot_1(x-y_0,x_0-x))\times (\bot_2(x_0,x+y-x_0)\cap \top_2(x+y-y_0,y_0))
\end{align*}
of $\sum_{y\in \Z_m}(\top(x,y)+\bot(x,y))$. 
Since $\sum_{y\in \Z_m}\top_1(x_0-y,y-y_0)=G_1$,  the coefficient of $(a,b)$ is equal to $m+1$ or $1$ depending on whether $a\in \bot_{1}(x-y_0,x_0-x)$ or not. 

Next, we assume that $b\in C_s^{(2m,q^4)}\cap C_t^{(q^2+1,q^4)}$ with $s$ and $t$ odd. Then,  there exists $(x_0,y_0)\in \Z_m\times \Z_m$ uniquely such that  
$b\in \bot_2(x,y_0)\cap \top_{2}(x_0,y)$ for all $x,y\in \Z_m$. 
Then,  $(a,b)$ may appear in the partial summation 
\begin{align*}
\sum_{y\in \Z_m}(\top_1(x-y_0,x_0-x)+\bot_1(x_0-y,y-y_0))
\times (\bot_2(x+y-y_0,y_0)\cap \top_2(x_0,x+y-x_0))
\end{align*}
of $\sum_{y\in \Z_m}(\top(x,y)+\bot(x,y))$. 
Since $\sum_{y\in \Z_m}\bot_1(x_0-y,y-y_0)=G_1$,  
the coefficient of $(a,b)$ is equal to 
$m+1$ or $1$ depending on whether $a\in \top_{1}(x-y_0,x_0-x)$ or not. 

Finally, we assume that $b=0$. Then, we have $b\in \bot_2(x,0)\cap \top_{2}(x',0)$ for all $x,x'\in \Z_m$. 
Then, $(a,b)$ may appear in the partial summation  
\begin{align*}
\sum_{y\in \Z_m}(\top_1(x,y)+\bot_1(x,y))\times (\bot_2(x+y,0)\cap \top_2(x+y,0))
\end{align*}
of $\sum_{y\in \Z_m}(\top(x,y)+\bot(x,y))$. 
Since $\sum_{y\in \Z_m}(\top_1(x,y)+\bot_1(x,y))$ satisfies the condition~(4), the coefficient of  $(a,b)$ is equal to  $m+1$, $2$, or $1$ depending on whether $a\in S_{x}$, $a=0$, or not. 

This completes the proof of the lemma. \qed
\vspace{0.3cm}

In the following lemma, we assume that $\top_1(x,y)$ and $\bot_1(x,y)$, $x,y\in \Z_m$,
satisfy the conditions~(1), (2), (3) and (5). 
(The proof is very similar to that of 
Lemma~\ref{lem:41}.) 
\begin{lemma}\label{lem:42}
The subsets $\top(x,y)$ and $\bot(x,y)$, where $x,y\in \Z_m$, defined in \eqref{eq:aaa1} and \eqref{eq:aaa2} satisfy the condition~(5). 
\end{lemma}
\pro 
Let $(a,b)\in G_1\times G_2$, where $a$ is arbitrarily chosen. 

First, we assume that  $b\in C_s^{(2m,q^4)}\cap C_t^{(q^2+1,q^4)}$ with $s$ and $t$ even. Then, there exists $(x_0,y_0)\in \Z_m\times \Z_m$ uniquely  such that  
$b\in \bot_2(x_0,y)\cap \top_{2}(x,y_0)$ for all $x,y\in \Z_m$. 
We now compute the coefficient of  $(a,b)$ in $\sum_{x\in \Z_m}(\top(x,y_1)+\bot(x,y_2))$. The $(a,b)$ may appear in the partial summation 
\begin{align*}
&\sum_{x\in \Z_m}\top_1(x_0-y_1,y_1-y_0)\times (\bot_2(x_0,-x_0+x+y_1)\cap \top_2(x+y_1-y_0,y_0))\\
+&\sum_{x\in \Z_m}\bot_1(x-y_0,x_0-x)\times (\bot_2(x_0,-x_0+x+y_2)\cap \top_2(x+y_2-y_0,y_0))
\end{align*}
of $\sum_{x\in \Z_m}(\top(x,y_1)+\bot(x,y_2))$. 
Since $\sum_{x\in \Z_m}\bot_1(x-y_0,x_0-x)=G_1$,  the coefficient of $(a,b)$ is equal to $m+1$ or $1$ depending on whether $a\in \top_{1}(x_0-y_1,y_1-y_0)$ or not. 

Next, we assume that 
$b\in C_s^{(2m,q^4)}\cap C_t^{(q^2+1,q^4)}$ with $s$ and $t$ odd. Then,  there exists $(x_0,y_0)\in \Z_m\times \Z_m$ uniquely such that  
$b\in \bot_2(x,y_0)\cap \top_{2}(x_0,y)$ for all $x,y\in \Z_m$. 
Then,  $(a,b)$ may appear in the partial summation 
\begin{align*}
&\sum_{x\in \Z_m}\top_1(x-y_0,x_0-x)
\times (\bot_2(x-y_0+y_1,y_0)\cap \top_2(x_0,-x_0+x+y_1))\\
+&\sum_{x\in \Z_m}\bot_1(x_0-y_2,-y_0+y_2)
\times (\bot_2(x+y_2-y_0,y_0)\cap \top_2(x_0,-x_0+x+y_2))
\end{align*}
of $\sum_{x\in \Z_m}(\top(x,y_1)+\bot(x,y_2))$. 
Since $\sum_{x\in \Z_m}\top_1(x-y_0,x_0-x)=G_1$,  
the coefficient of $(a,b)$ is equal to 
$m+1$ or $1$ depending on whether $a\in \bot_{1}(x_0-y_2,-y_0+y_2)$ or not. 

Finally, we assume that $b=0$. Then, we have $b\in \bot_2(x,0)\cap \top_{2}(x',0)$ for all $x,x'\in \Z_m$. 
Then, $(a,b)$ may appear in the partial summation  
\begin{align*}
&\sum_{x\in \Z_m}\top_1(x,y_1)\times (\bot_2(x+y_1,0)\cap \top_2(x+y_1,0))\\
+&\sum_{x\in \Z_m}\bot_1(x,y_2)\times (\bot_2(x+y_2,0)\cap \top_2(x+y_2,0))
\end{align*}
of $\sum_{x\in \Z_m}(\top(x,y_1)+\bot(x,y_2))$. 
Since $\sum_{x\in \Z_m}(\top_1(x,y_1)+\bot_1(x,y_2))$ satisfies the condition~(5), the coefficient of  $(a,0)$ with $a\not=0$ is equal to  $m+1$ or $1$. On the other hand, the 
coefficient of  $(0,0)$ is equal to  $0$, $2m$, or $m$ depending on whether $y_1\not=0$ and $y_2\not=0$, $y_1=y_2=0$, or not. 

This completes the proof of the lemma. \qed
\vspace{0.3cm}

Before proving our main theorems, we need to mention 
 about van Dam's theorem on ``fusions'' of strongly regular graphs.  
\begin{theorem}\label{thm:van}{\em (\cite{dam})}
Let $\Gamma_i$, $i=1,2,\ldots,h$, be strongly regular graphs decomposing the 
complete graph $K_v$. If $\Gamma_i$ are all of Latin square type or all of 
negative Latin square type, any union of $\Gamma_i$'s is also strongly regular. 
\end{theorem}
We are now ready for proving our main theorems. \vspace{0.2cm}

{\bf Proof of Theorem~\ref{thm:main1}: }
By Lemma~\ref{lem:41} and Proposition~\ref{prop:1234} together 
with Proposition~\ref{prop:inex}, 
each $S_x$, $x\in \Z_m$, generates a strongly regular Cayley graph on $G_1\times G_2\times \cdots\times G_s$ with negative Latin square type parameters $(u^2,c(u+1),- u+c^2+3 c,c^2+ c)$, where  $u=q_1^2q_2^2\cdots q_s^2$ and $c=(u-1)/m$. 
By Propositions~\ref{prop:1234}, the sets $S_x$, $x\in \Z_m$, partition $G\setminus \{0_G\}$, i.e., the complete graph $K_v$ is decomposed by the strongly regular Cayley graphs with connection sets 
$S_x$, $x\in \Z_m$.  
Then, by Theorem~\ref{thm:van}, the conclusion follows.  
\qed
\vspace{0.3cm}

{\bf Proof of Theorem~\ref{thm:main2}: }
By Lemma~\ref{lem:42} and Proposition~\ref{prop:1235} together 
with Proposition~\ref{prop:inex}, each $T_{y_1,y_2}$, $y_1,y_2\in \Z_m$, generates a strongly regular Cayley graph on $G_1\times G_2\times \cdots\times G_s$ with 
Latin square type parameters $(u^2,c(u-1), u+c^2-3 c,c^2- c)$, where $u=q_1^2q_2^2\cdots q_s^2$ and 
\[
c=
\begin{cases}
(u-1)/m, & \mbox{ if $y_1\not=0$ and $y_2\not=0$,}\\
(u-1)/m+2, & \mbox{ if $y_1=y_2=0$,}\\
(u-1)/m+1, & \mbox{ otherwise.}
\end{cases}
\]
By Propositions~\ref{prop:1235}, the sets $T_{y_1,y_1+y_2}$, $y_1\in \Z_m$, partition $G\setminus \{0_G\}$, i.e.,  the complete graph $K_v$ is decomposed by the strongly regular Cayley graphs with connection sets 
$T_{y_1,y_1+y_2}$, $y_1\in \Z_m$.  
Then, by Theorem~\ref{thm:van}, the conclusion follows. 
\qed

\section{Concluding remarks}
In this paper, we gave constructions of strongly regular Cayley graphs 
on abelian groups not necessarily elementary abelian groups by extending the theory of building blocks. 
Note that building blocks have rich applications in Combinatorics as mentioned in Introduction. 
For example, as a recent result, 
a building block of special type 
 (including our new building blocks) 
gives rise to a $1\frac{1}{2}$-difference set~\cite{DO}. Furthermore, the 
existence of a $1\frac{1}{2}$-difference set yields a 
directed strongly regular graph using its antiflags~\cite{BOS}. 
Thus, Theorem~\ref{thm:1half} itself is also important.

Next, we mention about designs related to our strongly regular graphs. It is known that the decomposition of the complete graph by  strongly regular graphs with same parameters $(v,k,\lambda,\mu)$ gives rise to a $2$-$(v,k,k-1)$ design~\cite{IM10}. In particular,  if the strongly regular graphs are Cayley graphs with connection sets $D_i$, where $i=1,2,\ldots,\frac{v-1}{k}$, 
then $\{D_i:i=1,2,\ldots,\frac{v-1}{k}\}$ forms a $(v,k,k-1)$ difference family in  the ambient group. Hence, by Theorem~\ref{thm:main1}, we can claim that there exists a $(v,k,k-1)$ difference family in $G_1\times G_2\times \cdots\times G_s$, where $v=q_1^4q_2^4\cdots q_s^4$ and $k=(v-1)/m$. As far as the author knows, this existence 
result of difference families is also new~\cite{AB}. 

Finally, we close this paper with an open problem. In this paper, we used a system of 
building blocks obtained from cyclotomic classes in semi-primitive case as starters of our product construction. Such a system of 
building blocks is a generalization of a covering extended building set found by Xia~(1992). However, no construction is known for a system of $2m^2$ building blocks satisfying the conditions (1), (2), (3) and either (4) or (5)  
in $(\F_{q^4},+)$ in the case where $m>2$ and $2m$ does not divide $q+1$.  On the other hand,   Chen~\cite{C97} found covering extended building sets in $(\F_{q}^4,+)$ for $q\equiv 1\,(\mod{4})$, where $2m$ does not divide  $q+1$ with $m=2$. Thus, we pose the following problem. 
\begin{problem}
Corresponding to Chen's result~\cite{C97}, 
find a construction of systems of $2m^2$ building blocks satisfying the conditions (1), (2), (3) and either (4)  or (5) in $(\F_{q}^4,+)$ 
in the case where $m>2$ and $2m$ does not divides $q+1$. 
\end{problem}



\begin{thebibliography}{99}

\bibitem{AB} R.~J.~R.~Abel, M.~Buratti, Difference families, in:  C.~J.~Colbourn,  
J.~H.~Dinitz (Eds.), {\em The CRC Handbook of Combinatorial Designs, 2nd edn},  392--409. Chapman~\&~Hall/CRC Press, Boca Raton, FL, (2006). 

\bibitem{BMW82}
L. D. Baumert, W. H. Mills, R. L. Ward, Uniform cyclotomy, {\it J. Number Theory} {\bf 14}, 67--82,  (1982).

\bibitem{bjl} T. Beth, D. Jungnickel, H. Lenz, {\it Design Theory}, Vol. I, 2nd edit., Cambridge University Press, 
Cambridge, (1999).

\bibitem{BH}
A. E. Brouwer, W. H. Haemers, {\it Spectra of Graphs}, Springer,  (2012).

\bibitem{BOS}
A. E. Brouwer, O. Olmez, S. Y. Song, Directed strongly regular graphs from $1\frac{1}{2}$-designs, {\it Europ. J. Comb.} {\bf 33}, 1174--1177, (2012) . 

\bibitem{BWX} A. E. Brouwer, R. M. Wilson, Q. Xiang, Cyclotomy and strongly regular graphs, 
{\it J. Algebraic Combin.} {\bf  10},  25--28, (1999).

\bibitem{C97}
Y. Q. Chen, On the existence of abelian Hadamard difference sets and a 
new family of difference sets, {\it Finite Fields Appl.} 
{\bf 3}, 234--256,  (1997). 


\bibitem{dam}
E. R. van Dam, Strongly regular decompositions of the complete graphs, 
{\it J. Algebraic Combin.} {\bf 17}, 181--201, (2003) 

\bibitem{DO}
J. A. Davis,  O. Olmez,  A framework for constructing partial geometric difference sets,  {\it Des. Codes, Cryptogr.} {\bf 86} 1367--1375, (2018).

\bibitem{DJ} J. A. Davis, J. Jedwab, A unifying construction for difference sets, {\it J. Combin. Theory, Ser. A,} {\bf 80}, 13--78, (1997). 



\bibitem{IM10}
T. Ikuta, A. Munemasa, Pseudocyclic association schemes and strongly regular graphs, {\em Europ. J. Combin.} {\bf 31}, 1513--1519,  (2010). 

\bibitem{LN97}
R. Lidl, H. Niederreiter, {\it Finite Fields}, Cambridge
Univ. Press, (1997).


\bibitem{Ma} 
S. L. Ma, A survey of partial difference sets, {\it Des. Codes Cryptogr.} 
{\bf 4}, 221--261, (1994).

\bibitem{MWX2019} 
K. Momihara, Q. Wang, Q. Xiang, Cyclotomy, difference sets, sequences with low
correlation, strongly regular graphs, and related
geometric substructures, in: K.-U. Schmidt, A. Winterhof (Eds.), {\em Combinatorics and Finite Fields. Difference Sets, Polynomials, Pseudorandomness and Applications}, 178--205 Radon Series on Computation and Applied Mathematics, {\bf 23}, De Gruyter (2019). 

\bibitem{P10}
J. Polhill, Paley partial difference sets in groups of 
order $n^4$ and $9n^4$ for any odd $n>1$,  
{\it J. Combin. Theory, Ser. A} {\bf 117}, 1027--1036,  (2010). 

\bibitem{T84}
R. J. Turyn, A special class of Williamson matrices and difference sets,  
{\it J. Combin. Theory, Ser. A} {\bf 36}, 111--115,  (1984). 

\bibitem{Wan2}
Z. Wang, New necessary conditions on (negative) Latin square
type partial difference sets in abelian groups, 
{\it J. Combin. Theory, Ser. A} {\bf 172}, 105208, (2020). 

\bibitem{Wan1}
Z. Wang, Paley type partial difference sets in abelian groups, 
{\it J. Combin. Designs} {\bf 28}, 149--152, (2020). 

\bibitem{WX97} 
R. M. Wilson, Q. Xiang, Constructions of Hadamard difference sets, 
 {\it J. Combin. Theory, Ser. A} {\bf 77}, 148--160,  (1997). 

\bibitem{X92}
M.-Y. Xia, Some infinite classes of special 
Williamson matrices and difference sets,  {\it J. Combin. Theory, Ser. A} {\bf 61}, 230-242,  (1992). 

\bibitem{XC96} 
Q. Xiang, Y. Q. Chen, On Xia's construction of Hadamard difference sets,  
 {\it Finite Fields Appl.}  
{\bf 2}, 87--95,  (1996). 

\end{thebibliography}
\end{document}